\documentclass[12pt]{article}
\usepackage{graphicx}%
\usepackage{multirow}%
\usepackage{amsmath,amssymb,amsfonts}%
\usepackage{amsthm}%
\usepackage{mathrsfs}%
\usepackage[title]{appendix}%
\usepackage{xcolor}%
\usepackage{textcomp}%
\usepackage{manyfoot}%
\usepackage{booktabs}%
\usepackage{algorithm}%
\usepackage{algorithmicx}%
\usepackage{algpseudocode}%
\usepackage{listings}%
\usepackage{dsfont}

\newcommand{\dm}{{\mathcal{D}}}
\newcommand{\R}{{\mathbb{R}}}
\newcommand{\N}{{\mathbb{N}}}
\newcommand{\p}{{\mathbb{P}}}
\newcommand{\D}{{\mathbb{D}}}

\newtheorem{theorem}{Theorem}
\newtheorem{defi}{Definition}

\newtheorem{cor}{Corollary}
\allowdisplaybreaks \allowdisplaybreaks[4]

\begin{document}

\title{Exponential bounds for the density of the law of the solution of a SDE with locally Lipschitz coefficients
\thanks{This work was supported by the Natural Sciences and Engineering Research Council of Canada under Grant DG-2018-04449.}
}

\author{Cristina Anton\\ 
Department of Mathematics and Statistics, MacEwan University, \\
103C, 10700-104 Ave., Edmonton, AB T5J 4S2, Canada \\
              Corresponding author, email: \texttt{popescuc@macewan.ca} 
}
\date{}
\maketitle

\begin{abstract}
Under the uniform  H\"{o}rmander's hypothesis we study smoothness and exponential bounds of the density of the law of the solution of a stochastic differential equation (SDE) with locally Lipschitz drift that satisfy a monotonicity condition. To obtain estimates for the Malliavin  covariance matrix and its inverse, we extend the approach in \cite{KusuokaStroock:1985} to SDEs with non-globally Lipschitz coefficients. As in \cite{ImkellerReisSalkeld:2019}, to avoid non-integrability problems we use results about Malliavin differentiability based on the concepts of  Ray Absolute Continuity and Stochastic Gate{\^a}ux differentiability.
\end{abstract}
\section{Introduction}\label{sec1}
We use Malliavin calculus to study smoothness and exponential bounds for the density  of the law of the solution of a stochastic differential equation (SDE) with a locally Lipschitz drift that satisfy a monotonicity condition. These exponential bounds are important, for example, to study the convergence rate of numerical schemes \cite{BallyTalay:1996} for approximating the solutions of the SDE. SDEs with non-globally Lipschitz coefficients appear in models for financial securities  and various models for dynamical systems such as stochastic Ginzburg–Landau equation and stochastic Duffing–van der Pol oscillator \cite{HutzenthalerJentzenKloeden:2011, KloedenPlaten:1992}.
  
We consider the SDE
\begin{equation}
dX(t)=b(X(t))dt+\sigma(X(t))dW(t),\quad
X(0)=x\in \mathbb{R}^d, t\in[0,T], T>0,\label{eq1}
\end{equation}
where $W(t)$ is an $m-$ dimensional Brownian motion defined on the filtered complete probability space $\left(\Omega, \mathcal{F},\{\mathcal{F}_t\}_{t\ge 0}, \mathbb{P}\right)$, and $b:\mathbb{R}^d\rightarrow \mathbb{R}^d $,  $\sigma:\mathbb{R}^d\rightarrow \mathbb{R}^{d\times m}$. 
We make the following assumptions for the coefficients  $b$ and $\sigma$:
\begin{itemize}
\item[\bf{C:}] $b$  has bounded partial derivatives of any order $k\ge 2$ and $\sigma$ has bounded partial derivatives of any order $k\ge 1$.
\item[\bf{M:}] There exist $L>0$, such that for any $x_1, x_2\in \mathbb{R}^d $ we have
\begin{equation}
<x_1-x_2,b(x_1)-b(x_2)>\le L|x_1-x_2|^2\label{c2}
\end{equation}
\item[\bf{P:}] There exists $L_1\ge 0$ and $N\ge 1$ such that for any $x_1, x_2\in\mathbb{R}^d$ we have
\begin{equation}
|b(x_1)-b(x_2)|^2\le L_1(1+|x_1|^{2N-2}+|x_2|^{2N-2})|x_1-x_2|^2\label{pol}
\end{equation}
\item[\bf{J:}]  There exist $L_3>0$ such that for any $y\in\mathbb{R}^d$ we have
\begin{equation}
y^\top \nabla_x by>-L_3|y|^2.
\end{equation}
\end{itemize}

Supposing that $b$ and $\sigma$ are globally Lipschitz, $C^\infty$,  all their derivatives have polynomial growth, and the H\"{o}rmander's hypothesis holds, in \cite{KusuokaStroock:1985} it is shown that  the strong solution of \eqref{eq1} is Malliavin differentiable of any order and it is nondegenerate at any fixed positive time. Furthermore, an estimate for the Malliavin covariance matrix (\cite[Theorem 2.17]{KusuokaStroock:1985}) is used to show that  the  law of the solution of the SDE is absolutely continuous with respect to the Lebesgue measure, its density is infinity differentiable and exponential bounds are proven under  the uniform H\"{o}rmander's hypothesis. 

There are several approaches to extend these results for SDE with non-globally Lipschitz coefficients. In \cite{DeMarco:2011}, assuming that the coefficients of the SDE are smooth and non-degenerate on an open domain $D$, estimations of the Fourier transform are used to show that the law of the solution has a smooth density and upper bounds for this density are given. In \cite{XieZhang:2016} the Sobolev regularity of strong solutions with respect to the initial value is established for SDEs with local Sobolev and super-linear growth coefficients. For SDEs driven by fractional Brownian motions, in \cite{BaudoinOuyangTindel:2014} it is shown that the density of the law of the solution is  smooth and admits an upper sub-Gaussian bound in the rough case.

For SDEs with random coefficients with drifts satisfying locally Lipschitz  and monotonicity conditions, in \cite{ImkellerReisSalkeld:2019} the concepts of Ray Absolute Continuity and Stochastic Gate{\^a}ux Differentiability are used to prove  Malliavin differentiability and absolute continuity of the solution's law. In \cite{Anton:2024arx} we extend this result  and  under assumptions  {\bf{C}}, {\bf{M}}, and {\bf{P}} we show Malliavin differentiability of any order. Here under assumptions  {\bf{C}}, {\bf{M}}, {\bf{P}} and {\bf{J}} we use the results in \cite{ImkellerReisSalkeld:2019} and \cite{Anton:2024arx} to get an estimate for the Malliavin covariance matrix similar with the one in (\cite[Theorem 2.17]{KusuokaStroock:1985}). If in addition the  uniform H\"{o}rmander's hypothesis holds, we  prove that the solution of the SDE is nondegenerate  and we obtain exponential bounds for the density of the  law of the solution of the SDE. 

The paper is organized as follows. In the next section we present some results regarding the  Malliavin differentiability of the solution of the SDE. Section \ref{sect3} includes estimates  for the Malliavin matrix.  Based on these estimates in section \ref{sect4} we prove  the exponential bounds for the density of the  law of the solution of the SDE \eqref{eq1}. 
\section{Notations and results about Malliavin differentiability}\label{sec2}
We denote by $\nabla f$ the gradient of a differentiable function $f:\mathbb{R}^d\rightarrow \mathbb{R}$, and for a vector valued function $v:\mathbb{R}^d\rightarrow \mathbb{R}^d$, let $\nabla v$ denote the matrix with components $\partial v_{i,j}(x)=\frac{\partial v^i(x)}{\partial x_j}$, $i,j=1,\ldots, d$.
For any multi-index $\alpha=(\alpha_1,\ldots,\alpha_d)$ with length $|\alpha|=\alpha_1+\cdots+\alpha_d$, let  $\partial^{\alpha}$ denote the partial derivative of order $|\alpha|$.
If $\phi$ is a smooth function, we denote by $\partial^\alpha_x\phi(t,x,y)$ the derivation with respect to the coordinates of $x$, where $t$, $y$ are fixed. 

For a vector $v\in\mathbb{R}^l$, we denote by $|v|:=\left(\sum_{i=1}^l v_i^2\right)^{\frac{1}{2}}$ the Euclidean norm, and if $A=(a_{ij})$ is an $l_1\times l_2$ matrix we denote by $\|A\|:=\left(\sum_{i=1}^{l_1}\sum_{j=1}^{l_2} a_{ij}^2\right)^{\frac{1}{2}}$ the Frobenius norm.
For two vectors $u,v\in\mathbb{R}^l$, we denote $<u,v>=\sum_{i=1}^l u_iv_i$, and for two $l_1\times l_2$ matrices  $A,B$,  $<A,B>=\sum_{i=1}^{l_1}\sum_{j=1}^{l_2} a_{ij}b_{ij}$ denotes the Frobenius inner product.

We consider the Banach space $\left(C([0,T]), \|\cdot\|_\infty\right)$ where
$C([0,T]):=\{ \phi:[0,T]\rightarrow \mathbb{R},~\phi \text{ uniformly continuous  and bounded}\}$, $\|\phi\|_\infty=\sup_{x\in [0,T]}|\phi(x)|$, and we denote $C_0([0,T]):=\{ \phi\in C([0,T]), \phi(0)=0\}$.

We say that a function $\phi:\mathbb{R}^l\rightarrow \mathbb{R}$ has polynomial growth, if there exist $C_f>0$ and $q_f\in \mathbb{N}_+$ such that for any $x\in \mathbb{R}^l$, we have
$|\phi(x)|\le C_f(1+|x|^{q_f})$. We define $\mathcal{C}^\infty_p(\mathbb{R}^l):=\{\phi\in \mathcal{C}^\infty(\mathbb{R}^l),~\phi \text{ and all its derivatives are functions with polynomial growth}\}$.

For any open set $E\subseteq \R^d$,  and $n\in \N$ we denote $C_b^n(E,\R^d)=\{ f\in C^n(E,\R^d), f$ and all its derivatives of order at most $n$ are bounded $\}$  with the norm $\|f\|_{C_b^n(E,\R^d)}=\max_{0\le i\le n}\sup_{x\in E}|\partial_x^i f(x)|$.

Let  $\left(\Omega, \mathcal{F},\{\mathcal{F}_t\}_{t\ge 0}, \mathbb{P}\right)$ be a filtered probability space. For any separable Banach space $(E,\|\cdot\|)$, we denote $L^p(\Omega;E)=\{X:\Omega\rightarrow E$ , $X$ is  $\mathcal{F}$ measurable and $\|X\|_p=E[\|X\|^p]^{1/p}<\infty  \}$. 
Let $L^\infty$ be the subset of bounded random variables with norm $\|X\|_{L^\infty}=ess \sup_{\omega\in \Omega}|X(\omega)| $.

Let $S^p([0,T],\mathbb{R}^d)=\{(Y_t)_{t\in [0,T]}$  stochastic processes, $Y_t\in\mathbb{R}^d$, that are $\{\mathcal{F}_t\}_{t\in [0,T]}$,   adapted and $\|Y\|_{S^p}=E[|Y\|^p_\infty]^{1/p}=E[\sup_{t\in [0,T]}|Y(t)|^p]^{1/p}$ $<\infty\}$. Let $S^\infty([0,T],\mathbb{R}^d)=\cap_{p\ge 1}S^p([0,T],\mathbb{R}^d)$.
\subsection{Malliavin calculus}
Let $\Omega=C_0\left([0,T],\mathbb{R}^m\right)=\{ \omega:[0,T]\rightarrow \mathbb{R}^m,\omega=(\omega_1,\ldots, \omega_m)^\top,$ $\omega_i\in C_0([0,T]), i=1,\ldots, m \} $ be the canonical Wiener space, and $W:=(W_t^1,\ldots, W_t^m)^\top_{t\in [0,T]}$ be the canonical Wiener process defined as $W_t^i(\omega):=\omega_t^i$ for any $\omega\in\Omega$, $i=1,\ldots, m$. We set $\mathcal{F}^0$ the natural filtration of $W$, $\mathbb{P}$ the Wiener measure, and $\mathcal{F}=\{\mathcal{F}_t\}_{t\in [0,T]}$ the usual augmentation (which is right-continuous and complete) of $\mathcal{F}^0$. In this setting $W$ is a standard Brownian motion. 

We denote $\mathcal{H}:=L^2([0,T];\mathbb{R}^m)=\{f:[0,T]\rightarrow \mathbb{R}^m$,  Borel measurable  and $\int_0^T|f(s)|^2 ds<\infty\}$ and the canonical inner product is
$$
<f,g>_{\mathcal{H}}:=\int_0^T <f(s),g(s)>ds=\sum_{i=1}^m\int_0^Tf^i(s)g^i(s)ds,~f, g\in \mathcal{H}.
$$
Let $H$ be the Cameron-Martin space:
\begin{align*}
H:=&\biggl\{h:[0,T]\rightarrow \mathbb{R}^m, h\in\Omega,\text{there exists } \dot{h}\in \mathcal{H}, \text{ such that } h(t)=\int_0^t \dot{h}(s)ds, t\in[0,T]\biggl\}
\end{align*}
For $h\in H$ we denote $\dot{h}$ a version of its Radon-Nykodym density with respect to the Lebesgue measure. For any Hilbert space $K$ we define $L^p(K)=\{f:\Omega\rightarrow K,$ $f$ is $\mathcal{F}_T$ measurable and $\|f\|^p_{L^p(K)}:=(E[\|f\|_K^p)^{1/p}<\infty\}$.  Let
$$
W(h):=\int_0^T \dot{h}_s dW_s:=\sum_{i=1}^m \int_0^T \dot{h}_s^idW_s^i,~h\in H.
$$
 Following \cite{MastroliaPossamaiRveillac:2017} we set
\begin{align*}
\mathcal{S}:=&\biggl\{F:\Omega\rightarrow \mathbb{R}, F=f(W(h_1), \ldots,W(h_n)), f\in C_b^\infty(\mathbb{R}^n), \\
&h_i\in H, i=1,\ldots, n, \text{ for some } n\in\mathbb{N}, n\ge 1\biggl\} 
\end{align*}
For any $F\in \mathcal{S}$ we define the Malliavin derivative $\dm F:\Omega\rightarrow H$ by
$$
\dm F:=\sum_{i=1}^n\partial_i f(W(h_1), \ldots,W(h_n))h_i.
$$
We identify $\dm F$ with the stochastic process $\{\dm_tF\}_{t\in[0,T]}$, where $\dm_tF\in \mathbb{R}^m$ and
$$
\dm_tF(\omega)=\sum_{i=1}^n\partial_i f(W(h_1)(\omega), \ldots,W(h_n)(\omega))h_i(t), ~(t,\omega)\in [0,T]\times \Omega.
$$
$\mathcal{D}_t^jF$ will denote the $j$th component of $\mathcal{D}_tF$. 
We denote $\mathbb{D}^{1,p}$, $p\ge 1$ the closure of $\mathcal{S}$ with respect to the semi-norm
$$
\|F\|_{1,p}:=\left(E[|F|^p]+E[\|\dm F\|_H^p\right)^{1/p},
$$
and we set $\mathbb{D}^{1,\infty}=\cap_{p\ge 2}\mathbb{D}^{1,p}$.

The $k$th order Malliavin derivative $\mathcal{D}^kF:\Omega\rightarrow H^k$ is defined iteratively and its components are
$\mathcal{D}_{t_1,\ldots, t_k}^{j_1,\ldots, j_k}F:=\mathcal{D}_{t_k}^{j_k}\ldots \mathcal{D}_{t_1}^{j_1}F$,  with $(t_1,\ldots, t_k)^\top\in [0,T]^k$, $j_1, \ldots, j_k\in \{1, \ldots, m\}$. For the $N$th order Malliavin derivative, $\mathbb{{D}}^{N,p}$, $p\ge 1$ is the closure of $\mathcal {S}$ with the semi-norm
$$
\|F\|_{N,p}:=\left(E\left[|F|^p\right]+\sum_{i=1}^N E\left[\|\dm^i F\|_{H^i}^p\right]\right)^{1/p}=\left(\|F\|^p_{L^p(\Omega)}+\sum_{i=1}^N\|\mathcal{D}^iF\|^p_{L^p\left(\Omega;L^2([0,T]^i,\mathbb{R}^m)\right)}\right)^{\frac{1}{p}}.
$$
We set $\mathbb{{D}}^\infty=\cap_{p\ge 2}\cap_{i\ge 1}\mathbb{{D}}^{i,p}$

The definition of Malliavin derivative can be extended to mappings $G:\Omega\rightarrow E$, where $(E,\|\cdot\|_E)$ is a separable Banach space (\cite{ImkellerReisSalkeld:2019}). We consider the family
\begin{align*}
\mathcal{S}_E:=&\biggl\{G:\Omega\rightarrow E, G=\sum_{j=1}^kF_je_j, F_j\in S,~e_j\in E \text{ for some } k\in\mathbb{N}, k\ge 1\biggl\} 
\end{align*}
$\mathcal{S}_E$ is dense in $L^p(\mathcal{F};E;\mathbb{P})$ \cite{ImkellerReisSalkeld:2019}.
For any $G\in \mathcal{S}_E$ we define the Malliavin derivative $\dm G:\Omega\rightarrow H\otimes E$ by
$$
\dm G:=\sum_{j=1}^k\dm F_j\otimes e_j
$$
We denote $\mathbb{D}^{1,p}(E)$, $p\ge 1$ the closure of $\mathcal{S}_E$ with respect to the semi-norm
$$
\|G\|_{1,p,E}:=\left(E[\|G\|_E^p]+E[\|\|\dm G\|_H\|_E^p\right)^{1/p},
$$

Next, following \cite{MastroliaPossamaiRveillac:2017} we will present a different characterization of the spaces $\mathbb{D}^{1,p}$. We start with some definitions from \cite{ImkellerReisSalkeld:2019}. Let $E$ be a separable Banach space and $L(H,E)$ be the space of all bounded linear operators $V:H\rightarrow E$.
\begin{defi}
A measurable map $f:\Omega\rightarrow E$ is said to be Ray Absolutely Continuous if for any $h\in H$ there exists a measurable mapping $\tilde{f}_h:\Omega\rightarrow E$ such that $\tilde{f}_h(\omega)=f(\omega)$, $\mathbb{P}$ a.s., and that for any $\omega\in \Omega$, $t\rightarrow \tilde{f}_h(\omega+th) $ is absolutely continuous on any compact subset of $\mathbb{R}$.
\end{defi}
\begin{defi}
A measurable map $f:\Omega\rightarrow E$ is said to be Stochastically Gate{\^a}ux differentiable if there exists a measurable mapping $F:\Omega\rightarrow L(H,E)$ such that for any $h\in H$,
$$
\frac{f(\omega+\epsilon h)-f(\omega)}{\epsilon}\overset{\mathbb{P}}{\rightarrow}F(\omega)[h]~\text{as }\epsilon\rightarrow 0.
$$
\end{defi}
\begin{defi}
Let $p>1$. A measurable map $f\in L^p(\Omega;E)$ is said to be Strong Stochastically Gate{\^a}ux differentiable if there exists a measurable mapping $F:\Omega\rightarrow L(H,E)$ such that for any $h\in H$,
$$
\underset{\epsilon\rightarrow 0}{\lim}E\left[\left\|\frac{f(\omega+\epsilon h)-f(\omega)}{\epsilon}-F(\omega)[h]\right\|\right]=0.
$$
\end{defi}
\begin{theorem}\cite[Theorem 3.10]{ImkellerReisSalkeld:2019}\\
\label{thg1}
Let $p>1$. The space $\mathbb{D}^{1,p}(E)$ is equivalent with the space of all random variables $f:\Omega\rightarrow E$ such that $f\in L^p(\Omega;E)$ is Ray Absolutely Continuous, Stochastically Gate{\^a}ux differentiable, and the Stochastic Gate{\^a}ux derivative $F:\Omega\rightarrow L(H,E)$ is $F\in L^p(\Omega;L(H,E)) $.
\end{theorem}
\begin{theorem}\cite[Theorem 3.13]{ImkellerReisSalkeld:2019}\\
\label{th3.13}
Let $p>1$. The space $\mathbb{D}^{1,p}(E)$ is equivalent with the space of all random variables $f:\Omega\rightarrow E$ such that $f\in L^p(\Omega;E)$ is Strong Stochastically Gate{\^a}ux differentiable and have measurable mappings  $F\in L^p(\Omega;L(H,E)) $.
\end{theorem}

\subsection{The solution of the SDE}\label{subsec2}
From assumptions {\bf{C}}, {\bf{M}}, {\bf{P}} we get for any  $x,y\in\R ^d$
\begin{align}
&y^\top\nabla_xb(x)y\le L|y|^2,\label{c4}\\ 
&|\sigma(x)|^2\le C(1+|x|^{2}), \label{us4}\\
&|\nabla_x \sigma(x)|^2\le C,\label{us1}\\
&|b(x)|^2\le C(1+|x|^{2N}),\label{us3}\\
&|\nabla_x b(x)|^2\le C(1+|x|^2).\label{us2}
\end{align}

From assumption {\bf{C}}, (\ref{us4}) and Theorems 3.6 in \cite{Mao:2011} we know that there exists a unique global solution $X(t,0,x)$ of the SDE (\ref{eq1}).
Also the solution $\left(X(t,0,x)\right)_{t\ge 0}$ is $\{\mathcal{F}_t\}_{t\ge 0}$ adapted and we have
\begin{equation}
E\left(\int_0^T |X(t,0,x)|^2dt\right)<\infty.\label{mom1}
\end{equation} 
Moreover from Theorems 9.1 and 9.5 in \cite{Mao:2011} we know that $\left(X(t,0,x)\right)_{t\ge 0}$ is a time homogeneous Markov process, and from Theorem 4.1 in \cite{Mao:2011} we know that  for any $p\ge  2$ there exists a constant $\alpha_p>0$ such that we have
\begin{equation}
E\left[|X(t,0,x)|^p\right]\le 2^{\frac{p-2}{2}}\left(1+|x|^p\right)e^{p\alpha_pt}, \quad t\in [0,T].\label{mom2}
\end{equation}

From Theorem 2.2 in \cite{ImkellerReisSalkeld:2019}, we know that the map $t\rightarrow X(t)(\omega)$, is $\mathbb{P}$ a.s. continuous, and  for any $p\ge 2$ we have $X\in S^p:=S^p([0,T],\mathbb{R}^d)$ and there exists $C>0$ depending on $p$, $b$ and $\sigma$ such that:
\begin{equation}
E\left[\sup_{t\in[0,T]}|X(t,0,x)|^p \right]< C(|x|^p+1).\label{mom3}
\end{equation}

From inequalities \eqref{us4}-\eqref{us2}, \eqref{mom3}, and assumption {\bf{C}}, for any multi-index $\alpha=(\alpha_1,\ldots,\alpha_d)$  and any $p\ge 2$ we get
\begin{align}
&E\left[\sup_{t\in[0,T]}|\partial_\alpha b(X(t))|^p \right]< \infty,~
E\left[\sup_{t\in[0,T]}|\partial_\alpha \sigma(X(t))|^p \right]< \infty.\label{mom31}
\end{align}

From Corollary 3.5 and Theorem 3.21 in \cite{ImkellerReisSalkeld:2019} we know that $X$ is Malliavin differentiable and  $\dm X\in S^p\left([0,T],L^2([0,T])\right)$ for any $p\ge 2$:
\begin{equation}
E\left[\left(\sup_{t\in[0,T]}\int_0^T|\dm_sX(t)|^2ds\right)^{p/2} \right]< \infty\label{malmom1}
\end{equation}
Thus from \eqref{mom3} and \eqref{malmom1}  we have $ X\in \mathbb{D}^{1,p}(S^p)$ for any $p\ge 2$, so $ X\in \mathbb{D}^{1,\infty}(S^p)=\cap_{p\ge 2}\mathbb{D}^{1,p}(S^p)$. Thus, similarly with the case in Theorem 2.2.1 in \cite{Nualart:2006} of globally Lipschitz coefficients, for any $t\in[0,T]$ we have $ X(t)\in \mathbb{D}^{1,\infty}$.

For any fixed $s\in [0,T]$ and $i=1,\ldots, m$, from \eqref{us4}, \eqref{c4}, \eqref{mom1}, and Theorem 2.5 in \cite{ImkellerReisSalkeld:2019} we have for any $p\ge 2$
\begin{align*}
&E\left[\underset{t\in [0,T]}{\sup} |\dm_s^iX(t)|^p\right]\le C E\left[|\sigma^i(X(s,0,x))|^p\right]\le C \left(1+E\left[|X(s))|^{p}\right]\right)
\end{align*}
This and \eqref{mom2} implies that for any $t\in [0,T]$ and any $p\ge 2$
\begin{align}
&\|X(t, 0,x)\|_{1,p}^p=E[|X(t, 0,x)|^p]+E\left[\left|\int_0^T|\dm_sX(t, 0,x)|^2 ds\right|^{p/2}\right]\le C_{1,p}(T)(1+|x|^p), \label{mom3-1}
\end{align}
where $C_{1,p}(T)>0$ depends on $p$, $T$, $b$ and $\sigma$.
In \cite{Anton:2024arx} we extend this result and show that under assumptions  {\bf{C}}, {\bf{M}}, and {\bf{P}}, $X^i(t)$ belongs to $\D ^\infty$ for all $t\in[0,T]$, and $i=1,\ldots, d$. Moreover, for any $t\in [0,T]$, $p\ge 2$, $k=1,2,\ldots, $ there exist $C_{k,p}(T), \beta_{k,p}>0$ depending on $p$, $k$, $T$, $b$, $\sigma$ such that
\begin{align}
\|X(t, 0,x)\|_{k,p}^p\le C_{k,p}(T)(1+|x|^{\beta_{k,p}})\label{dolfif}
\end{align}. 
\section{Results about the Malliavin matrix}
\label{sect3}
From Theorem 4.9 in \cite{ImkellerReisSalkeld:2019} we know that under assumptions {\bf{C}}, {\bf{M}}, {\bf{P}} the matrix valued SDE
\begin{equation}
J(t)=I_d+\int_0^t\nabla_x b\left(X(s,0,x)\right)J(s)ds+\int_0^t\nabla_x\sigma\left(X(s,0,x)\right)J(s)dW(s), \label{jaco}
\end{equation}
$t\in[0,T]$, has a unique solution $J\in S^p\left([0,T],\mathbb{R}^{d\times d}\right)$, $p\ge 2$, and for any $t\in[0,T]$ the map $x\rightarrow X(t,0,x)$ is differentiable $\mathbb{P}$ a.s. and as $\epsilon\rightarrow 0$
\begin{equation}
\frac{X(t, 0,x+\epsilon h)(\omega)-X(t,0,x)(\omega)}{\epsilon}\rightarrow hJ(t)(\omega)~\mathbb{P}~ a.s. . 
\end{equation}
From Theorem 2.5 and Proposition 4.13 in \cite{ImkellerReisSalkeld:2019} we know that under assumptions {\bf{C}}, {\bf{M}}, {\bf{P}}, {\bf{J}} the matrix valued SDE
\begin{align}
&K(t)=I_d-\int_0^tK(s)\left[\nabla_x b\left(X(s,0,x)\right)-<\nabla_x\sigma,\nabla_x\sigma>\left(X(s,0,x)\right)\right]ds\notag\\
&-\int_0^tK(s)\nabla_x\sigma\left(X(s,0,x\right)dW(s),\quad t\in[0,T], \label{jacin}
\end{align}
has a unique solution $K\in S^p\left([0,T],\mathbb{R}^{d\times d}\right)$, $p\ge 2$, and we have $K(t)J(t)=I_d$ for all $t\in[0,T]$ $\mathbb{P}$ a.s.. Consequently, the Jacobian matrix $J(t)$ is $\mathbb{P}$ a.s. invertible for any choice of $t\in[0,T]$, and $J(t)^{-1}=K(t)$ $\mathbb{P}$ a.s..

Let $J_s(t)=J(t)J(s)^{-1}$, $t>s$. Under  assumptions {\bf{C}}, {\bf{M}}, {\bf{P}}, {\bf{J}} from Proposition 5.1 in \cite{ImkellerReisSalkeld:2019} we know that we have
\begin{equation}
J_s(t)=I_d+\int_s^t\nabla_x b\left(X(r,0,x)\right)J_s(r)dr+\int_s^t\nabla_x\sigma\left(X(r,0,x)\right)J_s(r)dW(r), \label{jacofun}
\end{equation}
and the Malliavin derivative of $X$ can be expressed for $t>s$ as $\dm_sX(t,0,x)=J_s(t)\sigma(X(s,0,x))$. 
The Malliavin matrix $Q(t)$ is defined by
\begin{align}
Q(t,x)&:=\int_0^t\dm_sX(t,0,x)\dm_sX(t,0,x)^\top ds=J(t)C(t,x)J(t)^\top\label{matQ}\\
C(t,x)&:=\int_0^t J(s)^{-1}\sigma(X(s,0,x))\sigma(X(s,0,x))^\top(J(s)^{-1})^\top ds\label{matc}
\end{align}

The Lie bracket of the $C^1(\mathbb{R}^d,\mathbb{R}^d)$ vector fields $V=\sum_{i=1}^d V^i \frac{\partial}{\partial x_i}$, $U=\sum_{i=1}^d U^i \frac{\partial}{\partial x_i}$ is defined as
$[V,U](x)=\partial U(x)V(x)-\partial V(x) U(x)$, where $\partial U=(\partial_iU^j)_{i,j=1,\ldots d}$, $\partial V=(\partial_iV^j)_{i,j=1,\ldots d}$ are the Jacobian matrices of $U$ and $V$ respectively. Let us denote $\sigma^0=b-\frac{1}{2}\sum_{i=1}^{m}\sum_{j=1}^{d}\sigma_j^{i}\partial_j \sigma^i$ and let $\sigma^0$, $\ldots$, $\sigma^m$ be the corresponding vector fields 
$$
\sigma^0(x)=\sum_{i=1}^d \sigma^0_i(x)\frac{\partial}{\partial x_i},\quad \sigma^j(x)=\sum_{i=1}^d \sigma^j_{i}(x)\frac{\partial}{\partial x_i}, \quad j=1,\ldots, m.
$$
We construct by recurrence the sets $\Sigma_0=\{\sigma^j, j=1,\ldots,m \}$, $\Sigma_k=\{[\sigma^j, V],  j=0,\ldots,m, V\in\Sigma_{k-1}\},~ k\ge 1$, $\Sigma_\infty=\cup_{k=1}^\infty \Sigma_k$.
We denote by $\Sigma_k(x)$ the subset of $\mathbb{R}^m$ obtained by freezing the variable $x\in \mathbb{R}^d$ in the vector fields of $\Sigma_k$. For $x\in\mathbb{R}^d$ we consider the H\"{o}rmander's hypothesis:\\
{\bf{H(x):}} The vector space $Span\{\Sigma_\infty(x)\}=\mathbb{R}^d$.\\

As in \cite[Appendix]{KusuokaStroock:1985} let $\mathcal{A}=\{\emptyset\}\cup \bigcup_{i=1}^\infty (\{0,\ldots,m\})^i$. Given $\alpha=(\alpha_1,\ldots, \alpha_i)\in \mathcal{A}\backslash\{\emptyset\}$, we define $\alpha_{*}=\alpha_i$ and 
\begin{equation*}
\alpha^{'}=\begin{cases}
\emptyset,&\text{ if } i=1\\
(\alpha_1,\ldots, \alpha_{i-1}), &\text{ if }i\ge 2.
\end{cases}
\end{equation*}
Given $\alpha\in \mathcal{A}$ set
\begin{align*}
&|\alpha|=\begin{cases}
0 &\text{ if } \alpha=\emptyset\\
i&\text{ if } \alpha\in(\{0,\ldots, m\})^i,
\end{cases}\quad
&\|\alpha\|=\begin{cases}
0 &\text{ if } \alpha=\emptyset\\
|\alpha|+\text{card}\{j: \alpha_j=0\}&\text{ if } |\alpha|\ge 1.
\end{cases}
\end{align*}
We  define $T_{(\alpha)}$ and $I^{(\alpha)}(t)$ inductively on $|\alpha|$ by
\begin{align*}
&T_{(\alpha)}(V)=\begin{cases}
V&\text{ if } \alpha=\emptyset\\
[\sigma^{\alpha_*},T_{(\alpha^{'})}(V)]&\text{ if } \alpha\ne \emptyset,
\end{cases}, \quad V\in C^\infty(\R ^d, \R^d)\\
&I^{(\alpha)}(t)=\begin{cases}
1 &\text{ if } \alpha=\emptyset,\\
\int_0^t I^{(\alpha^{`})}(s)\circ dW^{\alpha^{*}}(s)&\text{ if } |\alpha|\ge 1,
\end{cases}
\end{align*}
where we consider $W^0(t)=t$, $t\in [0,T]$. 

Given $L\ge 1$ we define for any $x,\eta\in \R^d$
\begin{align*}
\mathcal{V}_L(x,\eta)&=\sum_{k=1}^m\sum_{\|\alpha\|\le L-1}<T_{(\alpha)}(\sigma^k)(x),\eta>^2, \quad
\mathcal{V}_L(x)&=\underset{|\eta|=1}{\inf}\mathcal{V}_L(x,\eta)\wedge 1.
\end{align*}
Let
\begin{align*}
U_L:=\{x\in\R^d, \mathcal{V}_L(x)>0\},\quad U=\bigcup_{L=1}^\infty U_L
\end{align*}
Notice that for $L\in \N^*=\{1,2, \ldots\}$ the hypothesis\\
H$_{\text{L}}$($x$): $Span\{\phi(x), \phi\in \cup_{i=1}^L \Sigma_i\}=\R^d$\\
is equivalent with $x\in U_L$. As in \cite{BallyTalay:1996} we consider the following assumption:\\
{\bf UH}: For some integer $L_0>0$, we have $C_{L_0}:=\inf_{x\in \mathbb{R}^d}V_{L_0}(x)>0$.\\
Notice that assumption {\bf UH} implies $U=\R^d$ and the H\"{o}rmander's hypothesis {\bf H(x)} is true for any $x\in \R^d$.

Suppose that  {\bf{H(x)}}, {\bf{C}}, {\bf{M}}, {\bf{P}}, {\bf{J}}  hold. Based on assumption {\bf{H(x)}}, \eqref{mom31}, the formulas \eqref{matQ}, \eqref{matc} for the Malliavin matrix $Q(t,x)$,  and proceeding as in the proof of Theorem 2.3.2 in  \cite{Nualart:2006} we can show that the Malliavin matrix $Q(t,x)$ is invertible a.s.. Thus, since from \eqref{mom3-1} we also know that $X(t, 0,x)\in \D ^{1,p}$ for any $p\ge 2$, this implies that  the law of $X(t,0, x)$ is absolutely continuous with respect to the Lebesgue measure (\cite[Theorem 2.2.1]{Nualart:2006}). Here we replace assumption {\bf{H(x)}} with assumption {\bf UH} and we obtain an exponential bound for the density of the law of $X(t,0, x)$ with respect to the Lebesgue measure.

Notice that we can write the equations \eqref{eq1}, \eqref{jaco}, \eqref{jacin} and \eqref{jacofun} in Stratonovich form as:
\begin{align}
&X(t)=x+\int_{0}^t \sigma^0(X(u)) du+\sum_{i=1}^m\int_{0}^t \sigma^i(X(u))\circ dW^i(u),\label{le1s}\\
&J(t)=I_d+\int_0^t\nabla_x \sigma^0\left(X(s)\right)J(s)ds+\sum_{i=1}^m\int_0^t\nabla_x\sigma^i\left(X(s)\right)J(s)\circ dW^i(s), \label{jacos}\\
&J^{-1}(t)=I_d-\int_0^t J^{-1}(s)\nabla_x \sigma^0\left(X(s)\right)ds-\sum_{i=1}^m\int_0^t J^{-1}(s)\nabla_x\sigma^i\left(X(s\right)) \circ dW^i(s), \label{jacins}\\
&J_s(t)=I_d+\int_s^t\nabla_x \sigma^0\left(X(r)\right)J_s(r)dr+\sum_{i=1}^m\int_s^t\nabla_x\sigma^i\left(X(r)\right)J_s(r)\circ dW^i(r). \label{jacofuns}
\end{align}
Given $V\in C^\infty(\R ^d, \R^d)$ we use Ito's formula and equations \eqref{le1s} and \eqref{jacins}   (\cite[equation 2.10]{KusuokaStroock:1985}, \cite[equation 2.63, pp 130]{Nualart:2006}) to obtain
\begin{align}
J^{-1}(t)V(X(t))&=V(x)+\int_{0}^t J^{-1}(s)[\sigma^0, V](X(s))ds\notag\\
&+\sum_{i=1}^m\int_0^t J^{-1}(s)[\sigma^i,V] (X(s))\circ dW^i(s)\label{expan1}
\end{align}

\begin{theorem}
\label{theorem212}
Suppose that assumptions {\bf{C}}, {\bf{M}}, {\bf{P}}, {\bf{J}}  hold. 
For any $L\ge 1$ and $0<\epsilon\le 1$ there exist $C_{L,\epsilon}, \lambda(L,\epsilon), \mu_{L,\epsilon}>0$ such that for all $x\in \R^d$ and $V\in C_p^\infty(\R^d,\R^d)$ we have
\begin{equation}
J^{-1}(t)V(X(t))=\sum_{\|\alpha\|\le L-1}T_{(\alpha)}(V)(x)I^{(\alpha)}(t)+R_L(t,x,V)\label{expJ3}
\end{equation}
for $K\ge K_0(L,x,V)> 1$, with
\begin{align}
&\sup_{0<t<1}\p\left(\frac{1}{t^L}\int_0^{t/K}|R_L(s,x,V)|^2ds\ge \frac{1}{K^{L+1-\epsilon}}\right)\le C_{L,\epsilon}\exp\left(\frac{-\lambda(L,\epsilon)K^{\mu_{L,\epsilon}}}{(1+M(x))^2}\right)
\end{align}
with 
\begin{align}
M(x)&=\max\{\|\sigma^i\|_{C_b^2(B(x,1),\R^d)},i=0,\ldots, m\}\vee \max\{\|T^{(\alpha)}(V)\|_{C_b^0(B(x,1),\R^d)}, |\alpha|\le L+1\}\notag
\end{align}
\end{theorem}
The proof is inlcuded in appendix \ref{secA1}.
\begin{theorem}
\label{th217}
Suppose that assumptions {\bf{C}}, {\bf{M}}, {\bf{P}}, {\bf{J}} hold.
For any $L\ge 1$  there exist $C(L)$, $\tilde{C}(L)>0$,  $\lambda(L)$, $\tilde{\lambda}(L)>0$, and $\mu_L$, $\tilde{\mu}_L\in (0,1]$, all of them independent of $\sigma^0, \ldots,\sigma^m$ such that for all $t\in (0,1]$ and all $K\ge 1$ we have
\begin{align}
&P\left(\frac{\tilde{\lambda}\left(\frac{t}{K^{1/(L+1)}},x\right)}{t^L}\le \frac{1}{K}\right)\le \tilde{C}(L)\exp\left(-\frac{\tilde{\lambda}_L\left(\mathcal{V}_L(x)^{L+2}K\right)^{\tilde{\mu}_L}}{(1+M(x))^2}\right)\label{lamba1}\\
&P\left(\frac{\lambda\left(\frac{t}{K^{1/(L+1)}},x\right)}{t^L}\le \frac{1}{K}\right)\le C(L)\exp\left(-\frac{\lambda_L\left(\mathcal{V}_L(x)^{L+2}K\right)^{\mu_L}}{(1+M(x))^2}\right),\label{lamba2}
\end{align}
where
\begin{align}
&\tilde{\lambda}(s,x)=\underset{|\eta|=1}{\inf}<\eta,C(s,x)\eta>,\quad \lambda(s,x)=\underset{|\eta|=1}{\inf}<\eta,Q(s,x)\eta>\label{lamb2}\\
&M(x)=\max\{\|T^{(\alpha)}(\sigma^k)\|_{C_b^2(B(x,1),\R^d)}, k=0,\ldots, m,|\alpha|\le L+1\}\label{Mx}
\end{align}
\end{theorem}
The proof is given in appendix \ref{secA2}.
\subsection{The Malliavin matrix under uniform H\"{o}rmander's hypothesis}
Let denote
$$
\Delta(t,x):=\det \left(Q(t,x)\right),\quad \tilde{\Delta}(t,x):=\det \left(C(t,x)\right)
$$
Equations \eqref{matc} and \eqref{lamb2} implies that $C(t,x)$ is positive semi-definite and both $C(t,x)$ and $\tilde{\lambda}(t,x)$ are non-decreasing with respect to $t$. Thus for $t\ge 1$
\begin{align}
Q(t,x)&=J(t)C(t,x)J^\top(t)\ge J(t)C(1,x)J^\top(t)\notag\\
&=J(t)J(1)^{-1}Q(1,x)J^\top(1)^{-1}J^\top(t)=J_1(t)Q(1,x)J_1^\top(t),\label{dety1}
\end{align}
where $J_1(t)=J(t)J(1)^{-1}$, $t\ge 1$ is the solution of \eqref{jacofun}.

From Theorem 2.5 and Proposition 5.1 in \cite{ImkellerReisSalkeld:2019} we know that under assumptions {\bf{C}}, {\bf{M}}, {\bf{P}} and {\bf{J}} the matrix valued SDE
\begin{align}
&G(t)=I_d-\int_1^tG(s)\left[\nabla_x b\left(X(s,0,x)\right)-<\nabla_x\sigma,\nabla_x\sigma>\left(X(s,0,x)\right)\right]ds\notag\\
&-\int_1^tG(s)\nabla_x\sigma\left(X(s,0,x\right)dW(s), \notag
\end{align}
$t\in[1,T]$, has a unique solution $G\in S^p\left([0,T],\mathbb{R}^{d\times d}\right)$, $p\ge 1$, and we have $G(t)J_1(t)=I_d$ for all $t\in[1,T]$ $\mathbb{P}$ a.s.. Consequently, the  matrix $J_1(t)$ is $\mathbb{P}$ a.s. invertible for any choice of $t\in[0,T]$, and $J_1(t)^{-1}=G(t)$ $\mathbb{P}$ a.s.. Moreover, for any $p\ge 1$ there exist $A_p$, $B_p>0$ depending on $C$ in \eqref{us1} such that for any $t\ge 1$ we have 
\begin{align}
E\left[\underset{u\in [1,t]}{\sup}\|J_1(u)^{-1}\|^{2p}\right]\le A_p\exp\left(B_p(t-1)\right)\label{sara1}
\end{align}
From \eqref{dety1} we get $Q(t,x)-J_1(t)Q(1,x)J_1^\top(t)\ge 0$ and since all the eigenvalues of a positive semi-definite matrix are non-negative we have $\det(Q(t,x))\ge \det(J_1(t))^2\det(Q(1,x))$. Since the Frobenius norm of a matrix $\|A\|\ge \sqrt{\lambda_{max}(AA^\top)}$ where $\lambda_{max}(AA^\top)$ is the largest eigenvalue of $AA^\top$, because $J(t)$  and $J(t)^{-1}$ are symmetric matrices we have 
\begin{align}
(\det(J_1(t)^{-1})^2&=\det(J_1(t)^{-1}(J_1(t)^{-1})^\top))\le\left( \lambda_{max}(J_1(t)^{-1}(J_1(t)^{-1})^\top)\right)^d \notag \\
&\le \|(J_1(t)^{-1}\|^{2d}.\label{detnorm}
\end{align}

\begin{theorem}
\label{thdelt}
Suppose that assumptions {\bf{C}}, {\bf{M}}, {\bf{P}}, and {\bf{J}} hold, and let $L\in \N^*$. For any $p\ge 1$ there exists $K(L,p)$, $\mu(L)>0$ also depending on the coefficents $b$ and $\sigma$, such that we have for any $t\in (0,1]$ and any $x\in U_L$
\begin{align}
E\left[\left|\frac{1}{\Delta(t,x)}\right|^p\right]\le K(L,p)\frac{(1+|x|^2)^{p\mu(L)}}{(\mathcal{V}_L(x)^{1+2/L}t)^{pdL}}\label{eqdety2}
\end{align}
\end{theorem}
\begin{proof}
Similarly with \eqref{detnorm} we get $(\det(J(t)^{-1})^2\le \|(J(t)^{-1}\|^{2d}$. Moreover we know from \eqref{matc} that $C(t,x)$  is positive semi-definite, so that the smallest eigenvalue of  $C(t,x)$ is equal with $\tilde{\lambda}(t,x)$,  and we have $(\tilde{\lambda}(t,x))^d\le \det(C(t,x))$.
These and \eqref{matQ} imply
\begin{align}
\left|\frac{1}{\Delta(t,x)}\right|&=\left|\frac{1}{\det(Q(t,x)}\right|=\left|\frac{(\det(J(t)^{-1})^2}{\det(C(t,x))}\right|\le \frac{(\det(J(t)^{-1})^2}{\tilde{\lambda}(t,x)^d}\le \frac{\|(J(t)^{-1}\|^{2d}}{\tilde{\lambda}(t,x)^d}\label{idet1bb}
\end{align}
From Theorem 2.5 in \cite{ImkellerReisSalkeld:2019} applied to the matrix valued SDE \eqref{jacins} we know that under assumptions {\bf{C}}, {\bf{M}}, {\bf{P}}, and {\bf{J}},  for any $p\ge 1$ there exist $a_p$, $b_p>0$ depending on $C$ in \eqref{us1} such that we have for any $t\in (0,1]$
\begin{align*}
E\left[\underset{u\in [0,t]}{\sup}\|J(u)^{-1}\|^{2p}\right]\le a_p\exp\left(b_pt\right)
\end{align*}
Replacing in \eqref{idet1bb}  and using Cauchy-Schwarz inequality we get for any $p\ge 1$ there exists $M_p>0$ such that
\begin{align}
E\left[\left|\frac{1}{(\Delta(t,x))^p}\right|\right]\le E\left[\|(J(t)^{-1}\|^{4pd}\right]^{1/2}E\left[\frac{1}{\tilde{\lambda}(t,x)^{2pd}}\right]^{1/2}\le M_p^{2pd}E\left[\frac{1}{\tilde{\lambda}(t,x)^{2pd}}\right]^{1/2}\label{fifi}
\end{align}
Next, since $t\rightarrow \tilde{\lambda}(t,x)$ is non-decreasing, from \eqref{lamba1} we get for any $t\in(0,1]$ and $K\ge 1$
\begin{align*}
&\p\left(\frac{\tilde{\lambda}(t,x)}{t^L}\le \frac{1}{K}\right)\le \p\left(\frac{\tilde{\lambda}(t/K^{1/(L+1)},x)}{t^L}\le \frac{1}{K}\right)\le \tilde{C}(L)\exp\left(-\frac{\tilde{\lambda}_L\left(\mathcal{V}_L(x)^{L+2}K\right)^{\tilde{\mu}_L}}{(1+M(x))^2}\right),
\end{align*}
where $\tilde{C}(L), \tilde{\lambda}_L>0$, $\tilde{\mu}_L\in (0,1]$,  and $M(x)$ are as in Theorem \ref{th217}.
Using this we get
\begin{align*}
&E\left[\frac{1}{\tilde{\lambda}(t,x)^{2pd}}\right]=\int_0^\infty 2pd y^{2pd-1}\p\left(\tilde{\lambda}(t,x)^{-1}>y\right)dy=\int_0^\infty 2pd y^{2pd-1}\p\left(\frac{\tilde{\lambda}(t,x)}{t^L}<\frac{1}{yt^L}\right)dy\\
&= \int_0^{1/A} 2pd y^{2pd-1}\p\left(\frac{\tilde{\lambda}(t,x)}{t^L}<\frac{1}{yt^L}\right)dy+\int_{1/A}^\infty 2pd y^{2pd-1}\p\left(\frac{\tilde{\lambda}(t,x)}{t^L}<\frac{1}{yt^L}\right)dy\\
&\le 2pd\left(\frac{1}{A}\right)^{2pd}+\int_{1/A}^\infty 2pd y^{2pd-1}\tilde{C}(L)\exp\left(-\frac{\tilde{\lambda}_L\left(Ay\right)^{\tilde{\mu}_L}}{(1+M(x))^2}\right)dy\\
&=2pd\left(\frac{1}{A}\right)^{2pd}+2pd\frac{\tilde{C}(L)}{\tilde{\mu}_L}\left(\frac{1}{A}\right)^{2pd}\int_{1}^\infty  z^{(2pd)/\tilde{\mu}_L-1}\exp\left(-\frac{\tilde{\lambda}_Lz}{(1+M(x))^2}\right)dz\\
&\le 2pd\left(\frac{1}{A}\right)^{2pd}+2pd\frac{\tilde{C}(L)}{\tilde{\mu}_L}\left(\frac{1}{A}\right)^{2pd}\int_{1}^\infty  z^{k-1}\exp\left(-\frac{\tilde{\lambda}_Lz}{(1+M(x))^2}\right)dz,
\end{align*}
where
\begin{align*}
A:=\mathcal{V}_L(x)^{L+2}t^L\in (0,1],\quad k=\left\lfloor \frac{2pd}{\tilde{\mu}_L}\right\rfloor+1\ge 1.
\end{align*}
By repeatedly applying integration by parts we get
\begin{align*}
&E\left[\frac{1}{\tilde{\lambda}(t,x)^{2pd}}\right]\le C_1(L,p)\left(\frac{1}{A}\right)^{2pd}(1+M(x))^{4pd/\tilde{\mu}_L+2}
\end{align*}
From assumption {\bf{C}} and \eqref{us4}-\eqref{us2} we get
\begin{align*}
M(x)^2\le C_1(1+|x|^{2N}).
\end{align*}
This yields
\begin{align*}
&E\left[\frac{1}{\tilde{\lambda}(t,x)^{2pd}}\right]\le C_1(L,p)\frac{(1+|x|^2)^{4Npd/\tilde{\mu}_L+2N}}{\left(\mathcal{V}_L(x)^{1+2/L}t\right)^{2pdL}}\le  C_1(L,p)\frac{(1+|x|^2)^{2Npd(2/\tilde{\mu}_L+1)}}{\left(\mathcal{V}_L(x)^{1+2/L}t\right)^{2pdL}}
\end{align*}
Thus for $\tilde{\mu}_L\in (0,1]$ replacing in \eqref{fifi} we get \eqref{eqdety2}.
\end{proof}
\begin{cor}
\label{cordet}
Suppose that assumptions  {\bf{C}}, {\bf{M}}, {\bf{P}}, {\bf{J}}, and {\bf{UF}} hold.   For any $p\ge 1$ and any $x\in \R^d$ there exists $L\in \N^*$, $\mu_L>0$, and $K_{L,p}(t)$ such that we have 
\begin{align}
E\left[\left|\frac{1}{\Delta(t,x)}\right|^p\right]\le K_{L,p}(t)\frac{(1+|x|^2)^{p\mu_L}}{t^{dpL}}, \quad t\in (0,T] \label{eqdety3}
\end{align}
with $t\rightarrow K_{L,p}(t)$ non-decreasing.
\end{cor}
\begin{proof}
We choose any $p\ge 1$ and any $x\in \R^d$. Form assumption {\bf{UF}} there exists $L\in \N^*$ such that $x\in U_L$. If $t\in (0,1]$, we get \eqref{eqdety3} from Theorem  \ref{thdelt}.

For $1< t\le T$, from \eqref{dety1} we get
$$
\left|\frac{1}{\Delta(t,x)}\right|=\left|\frac{1}{\det(Q(t,x))}\right|\le\frac{1}{\det(J_1(t))^2\det(Q(1,x))}=\frac{\det(J_1(t)^{-1})^2}{\Delta(1,x)}\le \frac{\|(J_1(t)^{-1}\|^{2d}}{\Delta(1,x)}
$$
Using Chauchy-Shwarz inequality and \eqref{sara1},  there exist $a_p$, $b_p>0$ depending on $C$ in \eqref{us1} such that
\begin{align}
&E\left[\left|\frac{1}{\Delta(t,x)}\right|^p\right]\le E\left[\left|\frac{\|J_1(t)^{-1}\|^{2d}}{\Delta(1,x)}\right|^p\right]\le E\left[\|J_1(t)^{-1}\|^{4dp}\right]^{1/2}E\left[\left|\frac{1}{\Delta(1,x)}\right|^{2p}\right]^{1/2}\notag\\
&\le a_{p}\exp\left(b_{p}(t-1)\right)E\left[\left|\frac{1}{\Delta(1,x)}\right|^{2p}\right]^{1/2}=a_{p}t^{pdL}\exp\left(b_{p}(t-1)\right)\frac{1}{t^{pdL}}E\left[\left|\frac{1}{\Delta(1,x)}\right|^{2p}\right]^{1/2} \notag
\end{align}
This and Theorem \ref{thdelt} imply \eqref{eqdety3}.
\end{proof}
\section{Exponential bounds for the density of  $X(t)$ under a uniform H\"{o}rmander's hypothesis}
\label{sect4}
\begin{theorem}
Let $X$ be the solution of SDE \eqref{eq1} and suppose that the assumptions {\bf{C}}, {\bf{M}}, {\bf{P}}, {\bf{J}}, and {\bf{UH}} hold. Then for any $t\in (0,T]$ and  any $x\in \R^d$ the law of the random vector $X(t,0,x)$ is absolutely continuous with respect to the Lebesgue measure and the density $y\rightarrow p_t(x,y)$  is a $C^\infty$ function. Moreover, for any $x\in \R^d$  the following inequalities hold
\begin{align}
&p_t(x,y)\le \frac{K_0(T)(1+|x|^{Q_0})}{t^{q_0}}\exp\left(-C_0\frac{(|x-y|\wedge 1)^2}{t(1+|x|)^{2N}}\right),\label{pf1}\\ 
&\left|\partial_y^\alpha p_t(x,y)\right|\le \frac{K_\alpha(T)(1+|x|^{Q_\alpha})}{t^{q_\alpha}}\exp\left(-C_\alpha\frac{(|x-y|\wedge 1)^2}{t(1+|x|)^{2N}}\right),\label{pf2}
\end{align}
for any $t\in (0,T]$, $y\in \R^d$, $t\le (|y-x|\wedge 1)/(4M(x))$, where $N$ is as in \eqref{us3} and $M(x)=\underset{z\in B(x,1)}{\sup}\{\|\sigma(z)\|\vee|b(z)|\}$. Here the non-decreasing functions $K_0$, $K_\alpha$   and the positive real numbers $C_0$, $C_\alpha$,  $Q_0$, $Q_\alpha$ depend on $L\in \N^{*}$ such that $x\in U_L$ and on the coefficients $b$, $\sigma$. 
\end{theorem}
\begin{proof}
 From Corollary \ref{cordet} we get that $1/\Delta(t,x)\in \bigcap_{p\ge 1}L^p(\Omega,\R)$, and since $X(t,0,x)\in \D^\infty$, $X(t,0,x)$ is non-degenerate for any $t\in (0,T]$, $x\in \R^d$ (see \cite[Definition 2.1.1]{Nualart:2006}). Thus we have the integration by parts formulas \cite[Proposition 2.1.4]{Nualart:2006}: for any $\phi\in C_p^\infty(R^d)$, $G\in \D^\infty$ and any index $\alpha$, there exists $H_\alpha (X(t,0,x),G)\in \D^\infty$ such that
 $$
 E[\partial^\alpha \phi(X(t,0,x))G]=E[\phi((X(t,0,x))H_\alpha(X(t,0,x),G)].
 $$
 Moreover,  for any $1< p<q_1<\infty$ there exist  constants $C_{p,q}>0$, $\beta, \gamma>1$ and $k_1, k_2\in \N^*$ such that  
 \begin{align}
 \|H_\alpha (X(t,0,x),G)\|_p\le C_{p,q_1}\left\|\frac{1}{\Delta(t,x)}\right\|_\beta^{k_1}\|\dm X(t,0,x)\|_{|\alpha|,\gamma}^{k_2}\|G\|_{|\alpha|,q_1}\label{hal}.
 \end{align}
Using this it can be shown \cite[Proposition 2.1.5]{Nualart:2006} that the density $p_t$ belongs to the Schwarz space $\mathcal{S}(\R^d)=\{f:\R^d\rightarrow \R: f\in C^\infty(\R^d), \underset{x\in \R^d}{\sup}|x|^k|\partial^\alpha f(x)|<\infty$ for any $k\ge 1$, and any index $\alpha$\}. Moreover,  for any $y=(y_1,\ldots,y_d)\in \R^d$, such that $y_i>0\vee x_i$, $i=1,\ldots, d$ we have
 \begin{align*}
 p_t(x,y)=E[\mathds{1}_{I_y}(X(t,0,x))H_{(1,\ldots,1)}(X(t,0,x),1)],\quad I_y=\prod _{i=1}^d[y_i,\infty).
 \end{align*}
 
 Using Cauchy-Schwarz inequality, \eqref{hal}, \eqref{eqdety3} in Corollary \ref{cordet},  and \eqref{dolfif} we get that for any $x\in \R^d$, $t\in (0,T]$   
 \begin{align}
 &p_t(x,y)\le E[\mathds{1}_{I_y}(X(t,0,x))]^{1/2}E[|H_{(1,\ldots,d)}(X(t,0,x),1)|^2]^{1/2}\notag\\
 &\le ( \p(X(t,0,x)\in I_y))^{1/2}C\left\|\frac{1}{\Delta(t,x)}\right\|_\beta^{k_1}\|\dm X(t,0,x)\|_{d,\gamma}^{k_2}\notag\\
 &\le C K_{L}(t)\frac{(1+|x|^Q)}{t^{dLk_1}}( \p(X(t,0,x)\in I_y))^{1/2},\label{f9}
 \end{align}
where $L\in \N^{*}$ is as in Corollary \ref{cordet}, $K_{L}(\cdot)$ is non-decreasing, $C>0$,  $\beta, \gamma>1$, and $k_1, k_2\in \N^*$.

Let 
\begin{align*}
\tau&=\underset{s\in[0,T]}{\inf}\{|X(s,0,x)-x|\ge 1/2\}\wedge T,\quad
\xi=\underset{s\in[0,T]}{\inf}\{|X(s,0,x)-x|\ge 1\}\wedge T
\end{align*}
If $|y-x|<1$ we have 
 \begin{align*}
& \p(X(t,0,x)\in I_y))=\p\left(\bigcap_{i=1}^d \{X_i(t,x)\ge y_i\}\right)=\p\left(\bigcap_{i=1}^d \{X_i(t,x)-x_i\ge y_i-x_i\}\right)\\
 &\le \p\left(|X(t,0,x)-x|\ge |y-x|\right)\\
 &\le\p\left(\underset{0\le s\le t\wedge \xi}{\sup}\left|\int_{0}^s b(X(u,0,x)) du+\sum_{i=1}^m\int_{0}^s \sigma^i(X(u,0,x))dW^i(u)\right|\ge |y-x|\right)\\
 &\le \p\left(\underset{0\le s\le t\wedge \xi}{\sup}\left(\int_{0}^s |b(X(u,0,x))| du+\sum_{i=1}^m\int_{0}^s |\sigma^i(X(u,0,x))|dW^i(u)\right)\ge |y-x|\right)
 \\
 \end{align*}
 For $s\le t\wedge \xi$ we have $|X(u,0,x)-x|\le 1$ for any $u\in [0,s]$,  so $X(u,0,x)\in B(x,1)$  and $|b(X(u,0,x))|\le M(x)$, $\|\sigma(X(u,0,x))\|\le M(x)$ for any $u\in [0,s]$. This implies
 $$
 \int_{0}^s \left|b(X(u,0,x))\right| du\le tM(x)\le |y-x|/2,
 $$
 for any $t\le |y-x|/(2M(x))$.
 We also have 
 $$
 \int_{0}^s |\sigma^i(X(u,0,x))|^2du\le M(x)^2t, \quad i=1,\ldots, m. 
 $$
 Hence for any $t\le |y-x|/(2M(x))$ we get
 \begin{align*}
 &\p\left(\underset{0\le s\le t\wedge \xi}{\sup}\left(\int_{0}^s |b(X(u,0,x))| du+\sum_{i=1}^m\int_{0}^s |\sigma^i(X(u,0,x))|dW^i(u)\right)\ge |y-x|\right)\notag\\
 &\le 2d\exp\left(-\frac{|y-x|^2}{8dtM(x)^2}\right)
 \end{align*}
 Here  we have applied Lemma 8.5 in \cite[Chapter V,Section 8]{IkedaWatanabe:1981} for each component of $X$.

 Similarly, if $|y-x|>1$
 \begin{align*}
& \p(X(t,0,x)\in I_y))\le \p\left(|X(t,0,x)-x|\ge |y-x|\right)\\
 &\le \p\left(\tau\le t\right)=\p\left(\underset{0\le s\le t\wedge \xi}{\sup}|X(s,0,x)-x|\ge 1/2\right)\\\\
 &=\p\left(\underset{0\le s\le t\wedge \xi}{\sup}\left|\int_{0}^s b(X(u,0,x)) du+\sum_{i=1}^m\int_{0}^s \sigma^i(X(u,0,x))dW^i(u)\right|\ge 1/2\right)\\
 &\le \p\left(\underset{0\le s\le t\wedge \xi}{\sup}\left(\int_{0}^s |b(X(u,0,x))| du+\sum_{i=1}^m\int_{0}^s |\sigma^i(X(u,0,x))|dW^i(u)\right)\ge 1/2\right)
 \end{align*}
 For $s\le t\wedge \xi$ we have $|X(u,0,x)-x|\le 1$ for any $u\in [0,s]$,  so $X(u,0,x)\in B(x,1)$  and $|b(X(u,0,x))|\le M(x)$, $\|\sigma(X(u,0,x))\|\le M(x)$ for any $u\in [0,s]$. This implies
 $$
 \int_{0}^s \left|b(X(u,0,x))\right| du\le tM(x)\le 1/4,
 $$
 for any $t\le 1/(4M(x))$.
 We also have 
 $$
 \int_{0}^s |\sigma^i(X(u,0,x))|^2du\le M(x)^2t, \quad i=1,\ldots, m. 
 $$
 Hence for any $t\le 1/(4M(x))$ we get
 \begin{align*}
 &\p\left(\underset{0\le s\le t\wedge \xi}{\sup}\left(\int_{0}^s |b(X(u,0,x))| du+\sum_{i=1}^m\int_{0}^s |\sigma^i(X(u,0,x))|dW^i(u)\right)\ge 1/2\right)\notag\\
 &\le 2d\exp\left(-\frac{1}{16dtM(x)^2}\right)
 \end{align*}
 Here we have applied Lemma 8.5 in \cite[Chapter V,Section 8]{IkedaWatanabe:1981} for each component of $X$.
 
 Thus we obtain for $t\le (|y-x|\wedge 1)/(4M(x))$
 \begin{align}
 \p(X(t,0,x)\in I_y))\le 2d\exp\left(-\frac{(|y-x|^2\wedge 1)}{16dtM(x)^2}\right)\label{f7}
 \end{align}
 Notice that from \eqref{us4}, \eqref{us3} we get
 \begin{align*}
 &M(x)^2\le C\underset{z\in B(x,1)}{\sup}(1+|z|)^{2N}\le C\underset{z\in B(x,1)}{\sup}(1+|z-x|+|x|)^{2N}\\
 &\le C(2+|x|)^{2N}\le 2^{2N}C(1+|x|)^{2N}
 \end{align*}
 Replacing in \eqref{f7} we obtain for $t\le (|y-x|\wedge 1)/(4M(x))$
 \begin{align}
 \p(X(t,0,x)\in I_y))\le C_1\exp\left(-C_2\frac{(|y-x|^2\wedge 1)}{t(1+|x|)^{2N}}\right)\label{f8}
 \end{align}
 This inequality and \eqref{f9} imply \eqref{pf1} for $t\le (|y-x|\wedge 1)/(4M(x))$.
 
 Next we prove \eqref{pf2}. From \cite[Proposition 2.1.5]{Nualart:2006} we know that for any index $\alpha$, and for any $y=(y_1,\ldots,y_d)\in \R^d$, such that $y_i>0\vee x_i$, $i=1,\ldots, d$ we have
 \begin{align*}
 \partial_y^\alpha p_t(x,y)=(-1)^{|\alpha|}E\left[\mathds{1}_{I_y}(X(t,0,x))H_\alpha\left(X(t,0,x),H_{(1,\ldots,1)}(X(t,0,x),1)\right)\right],
 \end{align*}
 where $I_y=\prod _{i=1}^d[y_i,\infty)$.
 Using this formula the proof is similar with the proof of \eqref{pf1}.
 \end{proof}



\appendix
\section*{Appendix A:Proof of Theorem \ref{theorem212}}\label{secA1}
\begin{proof}
The proof is similar with the proof of \cite[Theorem 2.12]{KusuokaStroock:1985}.
By repeated application of \eqref{expan1} we see that we have
\begin{align}
&J^{-1}(t)V(X(t))=V(x)+\sum_{i=0}^m\int_0^t J^{-1}(s)[\sigma^i,V] (X(s))\circ dW^i(s)=V(x)\notag\\
&+\sum_{i=0}^m[\sigma^i,V] (x)\int_0^t 1 \circ dW^i(s)+\sum_{i,j=0}^m\int_0^t \int_0^s J^{-1}(s)[\sigma^j,[\sigma^i,V] ](X(s))\circ dW^j(u)\circ dW^i(s)\notag\\
&=V(x)+\sum_{i=0}^m[\sigma^i,V] (x)\int_0^t 1 \circ dW^i(s)+\sum_{i,j=0}^m[\sigma^j,[\sigma^i,V] ](x)\int_0^t \int_0^s 1 \circ dW^j(u)\circ dW^i(s)\notag\\
&+\sum_{i,j,k=0}^m\int_0^t \int_0^s\int_0^u J^{-1}(s)[\sigma^k,[\sigma^j,[\sigma^i,V] ]](X(s))\circ dW^k(v)\circ dW^j(u)\circ dW^i(s)\notag\\
&=\sum_{|\alpha|\le L-1}T_{(\alpha)}(V)(x)I^{(\alpha)}(t)+\sum_{|\alpha|=L}S^{(\alpha)}(t,Z_{(\alpha)})\label{expJ2}
\end{align}
where as in \cite[Proof of Theorem 2.12]{KusuokaStroock:1985} we set (see also \cite[equation 2.63, pp 130]{Nualart:2006}) 
\begin{align*}
&S^{(\alpha)}(t,Z_{(\alpha)})=\begin{cases}
Z_{(\alpha)}(t) &\text{ if } \alpha=\emptyset,\\
\int_0^t S^{(\alpha^{`})}(s,Z_{(\alpha)})\circ dW^{\alpha^{*}}(s)&\text{ if } |\alpha|\ge 1,
\end{cases}\\
&Z_{(\alpha)}(t)=J^{-1}(t)T_{(\alpha)}(V)(X(t))=T_{(\alpha)}(V)(x)+\sum_{i=0}^m\int_0^t J^{-1}(s)[\sigma^i,T_{(\alpha)}(V)] (X(s))\circ dW^i(s)\\
&=T_{(\alpha)}(V)(x)+\sum_{i=1}^m\int_1^t J^{-1}(s)[\sigma^i,T_{(\alpha)}(V)] (X(s))dW^i(s)\\
&+\int_0^t J^{-1}(s)\biggl([\sigma^0,T_{(\alpha)}(V)] +\frac{1}{2}\sum_{i=1}^m[\sigma^i,[\sigma^i, T_{(\alpha)}(V)]]\biggl)(X(s))ds
\end{align*}
Notice that $S^{(\alpha)}(t,1)=I^{(\alpha)}(t)$.

From \eqref{expJ2} we have expansion \eqref{expJ3} with 
\begin{align*}
R_L(t,x,V)=\sum_{|\alpha|=L}S^{(\alpha)}(t,Z_{(\alpha)})+\sum_{\underset{|\alpha|\le L-1}{\|\alpha\|\ge L}}T_{(\alpha)}(V)(x)I^{(\alpha)}(t)
\end{align*}
For $f\in C([0,1])$ we have for any $t\in (0,1]$ and $K\ge 1$
\begin{align}
&\int_0^{t/K}|f(s)|^2ds= \int_0^{t/K}\left(\frac{|f(s)|}{s^{L/2-\epsilon/4}}\right)^2 s^{L-\epsilon/2}ds
\le \sup_{0<s\le 1/K}\left(\frac{|f(s)|}{s^{L/2-\epsilon/4}}\right)^2\int_0^{t/K} s^{L-\epsilon/2}ds\notag\\
&=\frac{2}{2L-\epsilon+2}\sup_{0<s\le 1/K}\frac{|f(s)|^2}{s^{L-\epsilon/2}}\left(\frac{t}{K}\right)^{L-\epsilon/2+1}\le \frac{2}{2L-\epsilon}\sup_{0<s\le 1/K}\frac{|f(s)|^2}{s^{L-\epsilon/2}}\left(\frac{t}{K}\right)^{L-\epsilon/2+1}\notag
\end{align}
Hence for any $0<t\le 1$, $K\ge 1$ we have 
\begin{align}
&\p\left(\frac{1}{t^L}\int_0^{t/K}|R_L(s,x,V)|^2ds\ge \frac{1}{K^{L+1-\epsilon}}\right)\le \p\left( \frac{2}{2L-\epsilon}\sup_{0<s\le 1/K}\frac{|R_L(s,x,V)|^2}{s^{L-\epsilon/2}}t^{1-\epsilon/2}\ge K^{\epsilon/2}\right)\notag\\
&\le \p\left( \frac{2}{2L-\epsilon}\sup_{0<s\le 1/K}\frac{|R_L(s,x,V)|^2}{s^{L-\epsilon/2}}\ge K^{\epsilon/2}\right)\le \p\left( \frac{2}{2L-\epsilon}\left(\sup_{0<s\le 1/K}\frac{|R_L(s,x,V)|}{s^{L/2-\epsilon/4}}\right)^2\ge K^{\epsilon/2}\right)\notag\\
&=\p\left( \sup_{0<s\le 1/K}\frac{|R_L(s,x,V)|}{s^{L/2-\epsilon/4}}\ge \left(\frac{2L-\epsilon}{2}\right)^{1/2}K^{\epsilon/4}\right)\notag\\
&\le \sum_{|\alpha|=L}\p\left( \sup_{0<s\le 1/K}\frac{|S^{(\alpha)}(s,Z_{(\alpha)})|}{s^{L/2-\epsilon/4}}\ge K_1\right)+\sum_{\underset{|\alpha|\le L-1}{\|\alpha\|\ge L}}\p\left( \sup_{0<s\le 1}\frac{|T_{(\alpha)}(V)(x)I^{(\alpha)}(s)|}{s^{L/2-\epsilon/4}}\ge K_1\right),\notag
\end{align}
where $K_1=\left(\frac{2L-\epsilon}{2}\right)^{1/2}\frac{K^{\epsilon/4}}{N}$, with $N=$ card$\{\alpha\in \mathcal{A}, |\alpha|=L\}+$card$\{\alpha\in \mathcal{A}, |\alpha|\le L-1, \|\alpha\|\ge L\}$.

To handle the terms of the first sum, let denote
\begin{align*}
Y_i^{(\alpha)}(s)&=J^{-1}(s)[\sigma^i,T_{(\alpha)}(V)] (X(s)),\quad i=1,\ldots, m\\
Y_0^{(\alpha)}(s)&=J^{-1}(s)\biggl([\sigma^0,T_{(\alpha)}(V)] +\frac{1}{2}\sum_{i=1}^m[\sigma^i,[\sigma^i, T_{(\alpha)}(V)]]
\end{align*}
We get for any $\alpha\in \mathcal{A}, |\alpha|=L$
\begin{align*}
&\p\left( \sup_{0<s\le 1/K}\frac{|S^{(\alpha)}(s,Z_{(\alpha)})|}{s^{L/2-\epsilon/4}}\ge K_1\right)\le \p\left( \sup_{0<s\le 1/K}\frac{|S^{(\alpha)}(s,Z_{(\alpha)})|}{s^{\|\alpha\|/2-\epsilon/4}}\ge K_1\right)\\
&\le \p\left( \sup_{0<s\le 1/K}\frac{|S^{(\alpha)}(s,Z_{(\alpha)})|}{s^{\|\alpha\|/2-\epsilon/4}}\ge K_1, \sup_{0<s\le 1/K}|Z_{(\alpha)}(s)|\le K_2, \biggl(\sum_{i=1}^m\int_0^{1/K}|Y_i^{(\alpha)}(s)|^2ds\biggl)^{1/2}\le K_2\right)\\
&+\p\left(\sup_{0<s\le 1/K}|Z_{(\alpha)}(s)|\ge K_2\right)+\p\left(  \biggl(\sum_{i=1}^m\int_0^{1/K}|Y_i^{(\alpha)}(s)|^2ds\biggl)^{1/2}\ge K_2\right)\\
&=:P_1(\alpha)+P_2(\alpha)+P_3(\alpha),
\end{align*}
where $K_2=K_1^{2^{-\|\alpha\|}}\ge K_1^{2^{-2L}}$ for $K$ large enough because $L\le \|\alpha\|\le 2|\alpha|=2L$.
From  \cite[Theorem A.5]{KusuokaStroock:1985} we know that there exist $C(L,\epsilon)>0$, $\lambda(L,\epsilon)>0$  such that 
\begin{align*}
P_1(\alpha)\le C(L,\epsilon)\exp\left(-\lambda_1(L,\epsilon)K_2\right)\le C(L,\epsilon)\exp\left(-\lambda_1(L,\epsilon)K_1^{2^{-2L}}\right).
\end{align*}
Thus there exist $C_1(L,\epsilon)>0$, $\lambda_1(L,\epsilon)>0$, $\mu_1(L,\epsilon)>0$  such that
\begin{align*}
\sum_{|\alpha|=L}P_1(\alpha)\le C_1(L,\epsilon)\exp\left(-\lambda_1(L,\epsilon)K^{\mu_1(L,\epsilon)}\right).
\end{align*}

For $K>1$, let
\begin{align*}
&\tau=\inf\{s\ge 0:\|J^{-1}(s)-I_d\|\vee|X(s)-x|\ge 1/2\}\wedge T\\
&\xi=\inf\{s\ge 0:\|J^{-1}(s)-I_d\|\vee|X(s)-x|\ge 1\}\wedge T
\end{align*}
We denote $F(s)=(X(s),J^{-1}(s))^\top$ We have
\begin{align*}
&\p\left(\sup_{0<s\le 1/K}|Z_{(\alpha)}(s)-Z_{(\alpha)}(0)|\ge K_2\right)\\
&\le \p\left(\sup_{0<s\le 1/K}|Z_{(\alpha)}(s)-Z_{(\alpha)}(0)|\ge K_2,\tau>1/K\right)+\p\left(\tau\le 1/K\right)\\
&\le\p\left(\sup_{0<s\le 1\wedge \tau}|Z_{(\alpha)}(s)-Z_{(\alpha)}(0)|\ge K_2\right)+\p\left(\sup_{0<s\le 1/K\wedge \xi}|F(s)-F(0)|\ge 1/2\right)
\end{align*}
For $K$ large enough  we  apply Lemma 8.5 in \cite[Chapter V, Section 8]{IkedaWatanabe:1981} for each component of $Z_{(\alpha)}$ and $F$  and we get:
\begin{align*}
&P_2(\alpha)\le  C_1 \exp \left(-\lambda_1K^2/(1+M(x))^2\right)+C_2\exp\left(-\lambda_2K/(1+M(x))^2\right)
\end{align*}

Notice that  we have
\begin{align*}
&P_3(\alpha)\le \sum_{i=1}^m\p\left(\int_0^{1/K}|Y_i^{(\alpha)}(s)|^2ds\ge \frac{K_2^2}{m}\right)\le \sum_{i=1}^m\p\left(\sup_{0<s\le 1/K}|Y_i^{(\alpha)}(s)|^2\ge \frac{K_2^2}{m}\right)\\
&\le \sum_{i=1}^m\p\left(\left(\sup_{0<s\le 1/K}|Y_i^{(\alpha)}(s)|\right)^2\ge \frac{K_2^2}{m}\right)=\sum_{i=1}^m\p\left(\sup_{0<s\le 1/K}|Y_i^{(\alpha)}(s)|\ge \frac{K_2}{\sqrt{m}}\right)
\end{align*}
For $Y_i^{(\alpha)}$, $i=1,\ldots, m$ we have a SDE similar with the one for $Z_{(\alpha)}$ but with $[\sigma^i,T_{(\alpha)}(V)]$ replacing $T_{(\alpha)}(V)$, so we can treat $P_2(\alpha)$ and the terms of $P_3(\alpha)$ similarly.

Finally since $S^{(\alpha)}(t,1)=I^{(\alpha)}(t)$,  we get for any $\alpha\in \mathcal{A}, |\alpha|\le L-1, \|\alpha\|\ge L$
\begin{align}
&\p\left( \sup_{0<s\le 1}\frac{|T_{(\alpha)}(V)(x)I^{(\alpha)}(s)|}{s^{L/2-\epsilon/4}}\ge K_1\right)\le \p\left( \sup_{0<s\le 1}\frac{|I^{(\alpha)}(s)|}{s^{L/2-\epsilon/4}}\ge \frac{K_1}{M(x)^2}\right)\notag\\
&\le \p\left( \sup_{0<s\le 1}\frac{|S^{(\alpha)}(s,1)|}{s^{\|\alpha\|/2-\epsilon/4}}\ge \frac{K_1}{2M(x)^2}\right)\notag
\end{align}
If $K$ is large enough, using \cite[Theorem A.5]{KusuokaStroock:1985} we know that there exist $C(L,\epsilon)>0$, $\lambda(L,\epsilon)>0$ and $\mu_1(L,\epsilon)>0$ such that 
\begin{align*}
&\sum_{\underset{|\alpha|\le L-1}{\|\alpha\|\ge L}}\p\left( \sup_{0<s\le 1}\frac{|T_{(\alpha)}(V)(x)I^{(\alpha)}(s)|}{s^{L/2-\epsilon/4}}\ge K_1\right)\le C(L,\epsilon)\exp\left(-\frac{\lambda(L,\epsilon)K^{\mu_1(L,\epsilon)}}{(1+M(x))^2}\right)
\end{align*}
\end{proof}
\section*{Appendix B:Proof of Theorem \ref{th217}}\label{secA2}
\begin{proof}
The proof is similar with the proof of \cite[Theorem 2.17]{KusuokaStroock:1985}.

From \eqref{matc} and \eqref{lamb2} notice that $C(t,x)$ is positive semi-definite and both $C(t,x)$ and $\tilde{\lambda}(t,x)$ are non-decreasing with respect to $t$, so it is enough to prove \eqref{lamba1} for $K\ge 1$ large enough. Since for any $a,b\ge 0$ we have $(a+b)^2\ge a^2/2-b^2$, from \eqref{matc}  and \eqref{expJ3} we have for any $t\in (0,1]$, $K>1$ and $\eta\in \R^d$, $|\eta|=1$
\begin{align*}
&<\eta,C(t/K,x)\eta>=\sum_{k=1}^m\int_0^{t/K}< J(s)^{-1}\sigma^k(X(s)),\eta>^2 ds\\
&=\sum_{k=1}^m\int_0^{t/K}<\sum_{\|\alpha\|\le L-1}T_{(\alpha)}(\sigma^k)(x)I^{(\alpha)}(s)+R_L(s,x,\sigma^k),\eta>^2ds\\
&=\sum_{k=1}^m\int_0^{t/K}\left(<\sum_{\|\alpha\|\le L-1}T_{(\alpha)}(\sigma^k)(x)I^{(\alpha)}(s),\eta>+<R_L(s,x,\sigma^k),\eta>\right)^2ds\\
&\ge \frac{1}{2}\sum_{k=1}^m\int_0^{t/K}<\sum_{\|\alpha\|\le L-1}T_{(\alpha)}(\sigma^k)(x)I^{(\alpha)}(s),\eta>^2ds-\sum_{k=1}^m\int_0^{t/K}<R_L(s,x,\sigma^k),\eta>^2ds\\
&\ge \frac{1}{2}\sum_{k=1}^m\int_0^{t/K}\left(<\sum_{\|\alpha\|\le L-1}T_{(\alpha)}(\sigma^k)(x),\eta>I^{(\alpha)}(s)\right)^2ds-\sum_{k=1}^m\int_0^{t/K}|R_L(s,x,\sigma^k)|^2ds
\end{align*}
This implies
\begin{align}
&\underset{|\eta|=1}{\inf}\sum_{k=1}^m\int_0^{t/K}\left(<\sum_{\|\alpha\|\le L-1}T_{(\alpha)}(\sigma^k)(x),\eta>I^{(\alpha)}(s)\right)^2ds\notag\\
&\le 2\underset{|\eta|=1}{\inf}<\eta,C(t/K,x)\eta>+2\sum_{k=1}^m\int_0^{t/K}|R_L(s,x,\sigma^k)|^2ds\notag\\
&=2\tilde{\lambda}(t/K,x)+2\sum_{k=1}^m\int_0^{t/K}|R_L(s,x,\sigma^k)|^2ds\label{lambd1}
\end{align}
Let 
$$
M_k(x,\eta):=\sum_{\|\alpha\|\le L-1}<T_{(\alpha)}(\sigma^k)(x),\eta>^2
$$
Since for any $\eta\in \R^d$, $|\eta|=1$
\begin{align*}
&\mathcal{V}_L(x,\eta)=\sum_{k=1}^m M_k(x,\eta)\ge \mathcal{V}_L(x),\quad
\sum_{\|\alpha\|\le L-1}\frac{<T_{(\alpha)}(\sigma^k)(x),\eta>^2}{M_k(x,\eta)}=1
\end{align*}
we get
\begin{align*}
&\underset{|\eta|=1}{\inf}\sum_{k=1}^m\int_0^{t/K}\left(<\sum_{\|\alpha\|\le L-1}T_{(\alpha)}(\sigma^k)(x),\eta>I^{(\alpha)}(s)\right)^2ds\\
&=\underset{|\eta|=1}{\inf}\sum_{k=1, M_k(x,\eta)>0}^m M_k(x,\eta) \int_0^{t/K}\left(\sum_{\|\alpha\|\le L-1}\frac{<T_{(\alpha)}(\sigma^k)(x),\eta>}{\sqrt{M_k(x,\eta)}}  I^{(\alpha)}(s)\right)^2ds\\
&\ge \underset{|\eta|=1}{\inf}\sum_{k=1}^m M_k(x,\eta)\inf\left\{\int_0^{t/K}\left(\sum_{\|\alpha\|\le L-1}b_\alpha  I^{(\alpha)}(s)\right)^2ds: \sum_{\|\alpha\|\le L-1} b_\alpha^2=1\right\}\\
&=\inf\left\{\int_0^{t/K}\left(\sum_{\|\alpha\|\le L-1}b_\alpha  I^{(\alpha)}(s)\right)^2ds: \sum_{\|\alpha\|\le L-1} b_\alpha^2=1\right\}\underset{|\eta|=1}{\inf}\mathcal{V}_L(x,\eta)\\
&\ge \inf\left\{\int_0^{t/K}\left(\sum_{\|\alpha\|\le L-1}b_\alpha  I^{(\alpha)}(s)\right)^2ds: \sum_{\|\alpha\|\le L-1} b_\alpha^2=1\right\}\mathcal{V}_L(x)
\end{align*}
Here if $M_k(x,\eta)=0$ then $<T_{(\alpha)}(\sigma^k)(x),\eta>=0$ for any $\|\alpha\|\le L-1$, so the term corresponding to such $k$ is 0.

Thus, for any $L \ge 1$, \eqref{lambd1} yields for any $\epsilon\in (0,1)$, $t\in (0,1]$, and $K\ge 1$
\begin{align*}
&\p\left(\frac{\tilde{\lambda}(t/K,x)}{t^L}\le \frac{1}{K^{L+1-\epsilon}}\right)\le \p\biggl(\frac{\mathcal{V}_L(x)}{2t^L}\inf\left\{\int_0^{t/K}\left(\sum_{\|\alpha\|\le L-1}b_\alpha  I^{(\alpha)}(s)\right)^2ds: \sum_{\|\alpha\|\le L-1} b_\alpha^2=1\right\}\\
&-\frac{1}{t^L}\sum_{k=1}^m\int_0^{t/K}|R_L(s,x,\sigma^k)|^2ds\le \frac{1}{K^{L+1-\epsilon}}\biggl)\le \p\biggl(\frac{1}{t^L}\sum_{k=1}^m\int_0^{t/K}|R_L(s,x,\sigma^k)|^2ds\ge \frac{1}{K^{L+1-\epsilon}}\biggl)\\
&+\p\biggl(\frac{\mathcal{V}_L(x)}{2t^L}\inf\left\{\int_0^{t/K}\left(\sum_{\|\alpha\|\le L-1}b_\alpha  I^{(\alpha)}(s)\right)^2ds: \sum_{\|\alpha\|\le L-1} b_\alpha^2=1\right\}\le \frac{2}{K^{L+1-\epsilon}}\biggl)\\
&\le \sum_{k=1}^m\p\biggl(\frac{1}{t^L}\int_0^{t/K}|R_L(s,x,\sigma^k)|^2ds\ge \frac{1}{mK^{L+1-\epsilon}}\biggl)\\
&+\p\biggl(\left(\frac{K}{t}\right)^L\inf\left\{\int_0^{t/K}\left(\sum_{\|\alpha\|\le L-1}b_\alpha  I^{(\alpha)}(s)\right)^2ds: \sum_{\|\alpha\|\le L-1} b_\alpha^2=1\right\}\le \frac{4}{\mathcal{V}_L(x)K^{1-\epsilon}}\biggl)
\end{align*}
By \cite[Theorem A.6]{KusuokaStroock:1985} there exist $C_1(L,\epsilon), \mu_1(L,\epsilon)>0$ such that for all $t\in (0,1]$
\begin{align*}
&\p\biggl(\left(\frac{K}{t}\right)^L\inf\left\{\int_0^{t/K}\left(\sum_{\|\alpha\|\le L-1}b_\alpha  I^{(\alpha)}(s)\right)^2ds: \sum_{\|\alpha\|\le L-1} b_\alpha^2=1\right\}\le \frac{4}{\mathcal{V}_L(x)K^{1-\epsilon}}\biggl)\\
&\le C_1(L,\epsilon)\exp\left(-\left(\frac{\mathcal{V}_L(x)K^{1-\epsilon}}{4}\right)^{\mu_1(L,\epsilon)}\right)
\end{align*}
By Theorem \ref{theorem212} there exist $C_2(L,\epsilon), \lambda_2(L,\epsilon), \mu_2(L,\epsilon)>0$ such that for all $t\in (0,1]$ and $K\ge 1$ large enough
\begin{align*}
&\sum_{k=1}^m\p\biggl(\frac{1}{t^L}\sum_{k=1}^m\int_0^{t/K}|R_L(s,x,\sigma^k)|^2ds\ge \frac{1}{mK^{L+1-\epsilon}}\biggl)\\
&\le m C_2(L,\epsilon)\exp\left(-\frac{\lambda_2(L,\epsilon)\left(m^{1/(L+1-\epsilon)}K\right)^{\mu_2(L,\epsilon)}}{(1+M(x))^2}\right)
\end{align*}
Replacing $K$ by $K^{1/(L+1)}$ and then taking $\epsilon=1/(L+2)$ and using $\mathcal{V}_L(x))\le 1$ we get that there exist  $\tilde{C}(L)>0$, $\tilde{\lambda}(L)>0$, and  $\tilde{\mu}_L>0$, all of them independent of $\sigma^0, \ldots,\sigma^m$ such that for all $t\in (0,1]$ and any $K\ge 1$ large enough
\begin{align*}
&\p\left(\frac{\tilde{\lambda}(t/K^{1/(L+1)},x)}{t^L}\le \frac{1}{K}\right)\le\p\left(\frac{\tilde{\lambda}(t/K^{1/(L+1)},x)}{t^L}\le \frac{1}{K^{(L+1-\epsilon)/(L+1)}}\right)\\
&\le  C_1(L,\epsilon)\exp\left(-\left(\frac{\mathcal{V}_L(x)K^{(1-\epsilon)/(L+1)}}{4}\right)^{\mu_1(L,\epsilon)}\right)\\
&+m C_2(L,\epsilon)\exp\left(-\frac{\lambda_2(L,\epsilon)\left(m^{1/(L+1-\epsilon)}K^{1/(L+1)}\right)^{\mu_2(L,\epsilon)}}{(1+M(x))^2}\right)\\
&\le C_1(L)\exp\left(-\left(\frac{\mathcal{V}_L(x)K^{1/(L+2)}}{4}\right)^{\mu_1(L)}\right)\\
&+m C_2(L)\exp\left(-\frac{\lambda_2(L)\left(m^{(L+2)/(L^2+3L+1)}K^{1/(L+1)}\right)^{\mu_2(L)}}{(1+M(x))^2}\right)\\
&\le C_1(L)\exp\left(-\frac{1}{4^{\mu_1(L)}}\frac{\left((\mathcal{V}_L(x))^{L+2}K\right)^{\mu_1(L)/(L+2)}}{(1+M(x))^2}\right)\\
&+m C_2(L)\exp\left(-\frac{\lambda_3(L)\left((\mathcal{V}_L(x))^{L+2}K\right)^{\mu_2(L)/(L+1)}}{(1+M(x))^2}\right)\\
&\le \tilde{C}(L)\exp\left(-\frac{\tilde{\lambda}_L\left(\mathcal{V}_L(x)^{L+2}K\right)^{\tilde{\mu}_L}}{(1+M(x))^2}\right)
\end{align*}
For $K\ge 1$ large enough  we can consider $\tilde{\mu}_L\in (0,1]$.

Next, notice that for any $\eta\in\R^d$, $|\eta|=1$
\begin{align*}
&<\eta, Q(s, x)\eta>=<J(s)^\top\eta, C(s,x)J(s)^\top\eta>\\
&=|J(s)^\top\eta|^2<\frac{J(s)^\top\eta}{|J(s)^\top\eta|}, C(s,x)\frac{J(s)^\top\eta}{|J(s)^\top\eta|}\ge |J(s)^\top\eta|^2 \tilde{\lambda}(s, x)
\end{align*}
Let
\begin{align*}
&\tau=\inf\{s\ge 0:\|J(s)-I_d\|\ge 1/2\}\wedge T,\quad \tau_1=\inf\{s\ge 0:\|J(s)-I_d\|\ge 1\}\wedge T\\
&\xi=\inf\{s\ge 0:|X(s)-x|\ge 1/2\}\wedge T,\quad\xi_1=\inf\{s\ge 0:|X(s)-x|\ge 1\}\wedge T.
\end{align*}
Notice that for $s\le \tau$ we have 
\begin{align*}
&|J(s)^\top\eta|^2=|(J(s)-I_d)^\top\eta+\eta|^2\ge \frac{1}{2}|\eta|^2-|(J(s)-I_d)^\top\eta|^2\\
&\ge \frac{1}{2}-\|J(s)-I_d\|^2|\eta|^2\ge \frac{1}{2}-\frac{1}{4}=\frac{1}{4}
\end{align*}
Hence for $s\le \tau$ and for any $\eta\in\R^d$, $|\eta|=1$
\begin{align*}
&<\eta, Q(s, x)\eta>\ge \frac{1}{4} \tilde{\lambda}(s, x)
\end{align*}
Thus there exist  ${C}(L)>0$, ${\lambda}(L)>0$, and  ${\mu}_L>0$, all of them independent of $\sigma^0, \ldots,\sigma^m$ such that for all $t\in (0,1)$ and all $K\ge 1$ large enough we have
\begin{align*}
&\p\left(\frac{\lambda\left(\frac{t}{K^{1/(L+1)}},x\right)}{t^L}\le \frac{1}{K}\right)\le \p\left(\frac{\tilde{\lambda}\left(\frac{t}{K^{1/(L+1)}},x\right)}{t^L}\le \frac{4}{K}, \tau>\frac{1}{K^{1/(L+1)}}\right)\\
&+\p\left(\tau\le \frac{1}{K^{1/(L+1)}}, \xi>\frac{1}{K^{1/(L+1)}}\right)+\p\left( \xi\le \frac{1}{K^{1/(L+1)}}\right)\\
&\le \p\left(\frac{\tilde{\lambda}\left(\frac{t}{K^{1/(L+1)}},x\right)}{t^L}\le \frac{4}{K}\right)+\p\left(\sup_{0<s\le \tau_1\wedge\frac{1}{K^{1/(L+1)}}}\|J(s)-I_d\|\ge 1/2, \xi>\frac{1}{K^{1/(L+1)}}\right)\\
&+\p\left( \sup_{0<s\le \xi_1\wedge\frac{1}{K^{1/(L+1)}}}|X(s)-x|\ge 1/2\right)\\
&\le C_0(L)\exp\left(-\frac{\lambda_0(L)\left(\mathcal{V}_L(x)^{L+2}K\right)^{\mu_0(L)}}{(1+M(x))^2}\right)+C_1(L)\exp\left(-\frac{\lambda_1(L)K^{1/(L+1)}}{(1+M(x))^2}\right)\\
&+C_2(L)\exp\left(-\frac{\lambda_2(L)K^{1/(L+1)}}{(1+M(x))^2}\right)\\
&\le C_0(L)\exp\left(-\frac{\lambda_0(L)\left(\mathcal{V}_L(x)^{L+2}K\right)^{\mu_0(L)}}{(1+M(x))^2}\right)+C_1(L)\exp\left(-\frac{\lambda_1(L)\left(\mathcal{V}_L(x)^{L+2}K\right)^{1/(L+1)}}{(1+M(x))^2}\right)\\
&+C_2(L)\exp\left(-\frac{\lambda_2(L)\left(\mathcal{V}_L(x)^{L+2}K\right)^{1/(L+1)}}{(1+M(x))^2}\right)\\
&\le C(L)\exp\left(-\frac{\lambda_L\left(\mathcal{V}_L(x)^{L+2}K\right)^{\mu_L}}{(1+M(x))^2}\right)\\
\end{align*}
Here we have used $\mathcal{V}_L(x)\le 1$, and for the first probability in the third inequality we have used \eqref{lamba1}(we have $4/K$ instead of $1/K$ but we can ajust the constants from the beginning of the proof of \eqref{lamba1}). For the second and third probabilitites in the the third inequality  for $K>1$ large enough we have applied Lemma 8.5 in \cite[Chapter V, Section 8]{IkedaWatanabe:1981} for each component of $X$ and $J$. 
\end{proof}

\end{document}